\documentclass[hidelinks,onefignum,onetabnum]{siamart220329}

\usepackage{graphicx}
\usepackage{amsfonts,amsbsy,amscd,amsgen,dsfont,amsmath,amssymb,mathrsfs}
\usepackage{bm,multirow,cancel}
\usepackage{tabularx}
\usepackage{enumitem}
\usepackage{algorithm,algpseudocode}
\usepackage{textcomp}

\newcommand{\mmin}{\overline\min}
\newcommand{\id}{\mathrm d}
\newcommand{\bxi}{\pmb \xi}
\newcommand{\vc}{\mathbf}

\newcommand{\etab}{\pmb{\eta}}

\renewcommand{\tilde}{\widetilde}

\newcommand{\Null}{\mathcal{N}}

\newcommand{\fm}{\varphi}

\DeclareMathAlphabet\mathbfcal{OMS}{cmsy}{b}{n}
\DeclareMathOperator*{\argmin}{arg\,min}

\newcommand{\bphi}{\pmb\phi}

\newcommand{\R}{\mathcal{R}}

\newtheorem{thm}{Theorem}
\newtheorem{rem}{Remark}
\newtheorem{lem}{Lemma}

\newtheorem{defn}{Definition}
\newtheorem{ass}{Assumption}
\newtheorem{cor}{Corollary}

\headers{Sparse Discrete Empirical Interpolation}{M. Farazmand}

\title{Sparse Discrete Empirical Interpolation Method: State Estimation from Few Sensors\thanks{Submitted to the editors on Feb. 1, 2024.
		\funding{This work was supported by the National Science Foundation, the Algorithms for Threat Detection (ATD) program, through the award DMS-2220548.}}
	}

\author{Mohammad Farazmand
	\thanks{Department of Mathematics, North Carolina State University,
		2311 Stinson Drive, Raleigh, NC 27695-8205, USA (\email{farazmand@ncsu.edu}).}
	}

\usepackage{amsopn}

\begin{document}

\maketitle

\begin{abstract}
Discrete empirical interpolation method (DEIM) estimates a function from its incomplete pointwise measurements. 
Unfortunately, DEIM suffers large interpolation errors when few measurements are available. Here, we introduce Sparse DEIM (S-DEIM) for accurately estimating a function even when very few measurements are available. To this end, S-DEIM leverages a kernel vector which has been neglected in previous DEIM-based methods. We derive theoretical error estimates for S-DEIM, showing its relatively small error when an optimal kernel vector is used. When the function is generated by a continuous-time dynamical system, we propose a data assimilation algorithm which approximates the optimal kernel vector using sparse observational time series. We prove that, under certain conditions, data assimilated S-DEIM converges exponentially fast towards the true state. We demonstrate the efficacy of our method on two numerical examples.
\end{abstract}



\section{Introduction}\label{sec:intro}
Estimating the full state of a system from partial observations is needed in many fields such as control theory, fluid mechanics, imaging, meteorology, and oceanography~\cite{Manohar2018}. Methods that seek to accomplish this task go by various names, such as state estimation~\cite{chernousko1993}, flow reconstruction~\cite{Zaki2021b}, compressed sensing~\cite{Donoho2006}, and optimal experimental design~\cite{Pukelsheim2006}. 

Among existing methods, Discrete Empirical Interpolation Method (DEIM) stands out for its simplicity and interpretability~\cite{Sorensen2010}. DEIM, and its later variant Q-DEIM~\cite{Drmac2016}, use pointwise measurements of the state of the system to reconstruct the full state as a linear combination of prescribed basis functions (or modes). As we review in section~\ref{sec:relWork}, (Q-)DEIM first approximates the optimal sensor placement for gathering pointwise observations. Given this data, it then determines the optimal expansion coefficients in the prescribed basis in order to estimate the full state. Our contributions here concern the latter step; we use the Q-DEIM algorithm~\cite{Drmac2016} for sensor placement.

Numerical evidence in~\cite{Brunton2021} suggests that (Q-)DEIM performs best when the number of sensors are greater than or equal to the number of modes. However, in the regime where the number of sensors are smaller than the number of modes, (Q-)DEIM can incur large reconstruction errors. Unfortunately, one frequently encounters this regime in practice since
sensors can be quite expensive and furthermore their deployment in large numbers may face practical obstacles~\cite{Brunton2021}.
On the other hand, modes or basis functions are either prescribed analytically (e.g., polynomial interpolation) or are obtain numerically (e.g., via proper orthogonal decomposition). As such, the number of modes can be increased arbitrarily as needed.

This paper addresses the latter regime where the number of sensors is smaller than the number of modes. We first introduce \emph{Sparse DEIM} (S-DEIM) by making use of the kernel vector which has been neglected in previous DEIM-based methods. Our error estimates show that certain values of this kernel vector result in much more accurate state estimation compared to (Q-)DEIM. However, determining the optimal kernel vector from observations is non-trivial. In the case where the state is generated by a continuous-time dynamical system, we propose a data assimilation method which approximates the optimal kernel vector from observational time series. Therefore, our contributions can be summarized as follows:
\begin{enumerate}
	\item Sparse DEIM: We introduce S-DEIM for state reconstruction when the number of sensors is smaller than the number of modes. 
	\item Error estimates : We derive a number of error estimates for S-DEIM which quantify its approximation accuracy and the optimal choice of the kernel vector.
	\item Ruling out two-stage state estimation: We prove that a two-stage state estimation using S-DEIM, although seemingly promising, fails to improve the (Q-)DEIM reconstruction.
	\item Data assimilation with S-DEIM: In the special case, where the state is generated by a dynamical system, we propose a data assimilation method for estimating the optimal kernel vector in S-DEIM. Under certain assumptions, we prove that our method converges exponentially fast towards the true state.
\end{enumerate}

\subsection{Related work}\label{sec:relWork}
Due to its broad applications, a plethora of methods have been developed for state reconstruction from sparse measurements. 
Broadly speaking, these methods can be divided into two categories. One category of methods is deterministic and relies on solving a least squares problem; DEIM belongs to this category. The second category is based on Bayesian inference and yields a probabilistic reconstruction~\cite{Attia2018,Ji2008,Klishin2023}.
Furthermore, recent years have seen a surge of methods based on deep learning~\cite{Asch2022,Brunton2020_ARFM} which overcome some drawbacks~\cite{Rowley2022} of linear methods such as DEIM. In this section, we only review previous work on (Q)-DEIM and refer to~\cite{Callaham2019,Farazmand2023} for a broader review of other methods.

DEIM was originally proposed by Chaturantabut and Sorensen~\cite{Sorensen2010} for reduced-order modeling of dynamical systems generated by nonlinear differential equations. Later, it became clear that the same methodology can also be used for estimating the full state of the system from its sparse sensor measurements~\cite{Manohar2018}.

The original work of Chaturantabut and Sorensen~\cite{Sorensen2010} also proposed a method for optimal sensor placement. However, this method was not invariant under reordering of the basis functions. Drma\v c and Gugercin~\cite{Drmac2016} proposed Q-DEIM to address this issue. Q-DEIM uses QR factorization with column pivoting to approximate the optimal placement of the sensors. Furthermore, the theoretical upper bound for Q-DEIM error is sharper than the one derived in~\cite{Sorensen2010}. 
There are several other generalizations of DEIM, including its application to weighted inner product spaces~\cite{Drmac2018}, tensor-based DEIM~\cite{Farazmand2023,Kirsten2022}, DEIM for classification~\cite{Hendryx2021}, and its generalization for state reconstruction on manifolds~\cite{Rowley2019} as opposed to a linear subspace.

(Q-)DEIM is well-known to be sensitive to measurement error. Peherstorfer et al.~\cite{Drmac2020} carry out a probabilistic analysis to quantify the reconstruction error of DEIM when the observations are noisy. A similar analysis, in the context of D-optimal design, is carried out in~\cite{Karnik2024,Manohar2018}. Callaham et al.~\cite{Callaham2019} propose several methods which regularize the underlying DEIM optimization problem in order to decrease its sensitivity to noise. Unlike DEIM, their modified optimization problems do not necessarily have a closed-form solution. As a result, their computational cost can be significantly higher than that of DEIM. 

DEIM was originally developed under the assumption that the number of sensors $n$ is equal to the number of modes $m$ used in the linear expansion~\cite{Sorensen2010}. However, DEIM can easily be generalized to the case where the number of sensor and modes are different ($n\neq m$). Numerical evidence suggests that  the (Q-)DEIM reconstruction is most accurate when the number of sensors is greater than or equal to the number of modes  ($n\geq m$)~\cite{Brunton2021}. This regression regime has been studied extensively.  In fact, when $n>m$, DEIM coincides with gappy proper orthogonal decomposition (POD) which was discovered much earlier~\cite{Astrid2008,Willcox2004,Sirovich1995,Willcox2006}. Peherstorfer et al.~\cite{Drmac2020} showed that gappy POD is more robust to observational noise compared to DEIM with $n=m$. Recently, Lauzon et al.~\cite{Shin2022} proposed a point selection algorithm for optimal placement of additional sensors in gappy POD (also see~\cite{Shin2016}).

The opposite regime, where fewer sensors are available compared to modes ($n<m$), has received relatively little attention. As we mentioned earlier, this regime arises frequently in practice since procurement and deployment of sensors can be quite expensive. The present work develops the method of Sparse DEIM to address this regime: fewer sensors than modes.

\subsection{Outline}\label{sec:outline}
The remainder of this paper is organized as follows. In section~\ref{sec:prelim}, we discuss the problem set-up and review some relevant aspects of (Q-)DEIM. In section~\ref{sec:S-DEIM}, we present our S-DEIM method and its error analysis. In section~\ref{sec:DAS-DEIM}, we present a computable method for S-DEIM in the case where the states are generated by a continuous-time dynamical system. Section~\ref{sec:numerics} contains our numerical examples. Our concluding remarks are presented in section~\ref{sec:conc}.

\section{Preliminaries and set-up}\label{sec:prelim}
We consider bounded functions $u: \Omega \to \mathbb R, \vc x\mapsto u(\vc x)$, where $\Omega\subset \mathbb R^d$ is the spatial domain over which the function $u$ is defined. 
We assume that the function $u$ belongs to a Banach space $\mathcal U$. Consider a complete basis $\{\phi_i\}_{i=1}^\infty$ in $\mathcal U$, so that
\begin{equation}
\lim_{m\to\infty} \sum_{i=1}^m c_i\phi_i(\vc x) = u(\vc x),
\label{eq:inf_series}
\end{equation}
for appropriate coefficients $c_i\in\mathbb R$, where the convergence is uniform in the norm $\|\cdot\|_{\mathcal U}$ associated with the Banach space $\mathcal U$.

Consider a finite truncation of the series in~\eqref{eq:inf_series}, 
\begin{equation}
\tilde u(\vc x) = \sum_{i=1}^m c_i\phi_i(\vc x),
\label{eq:fin_series}
\end{equation}
for some $m>0$.
We would like to find the coefficients $c_i$ such that $\tilde u \simeq u$, where the approximation is understood in an appropriate sense. In order to find the coefficients $c_i$ some additional information about the function $u$ is needed.
We assume that the function values are measured at a sparse set of distinct points $\mathcal S:=\{\vc s_i\}_{i=1}^{n}\subset\Omega$. In practice, these measurements are obtained by $n$ sensors located at the points $\vc s_i$. Given the measurements $u(\vc s_i)$, we would like to reconstruct the function $u$ over the entire domain $\Omega$.

Therefore, the problem statement can be informally summarized as follows: Given the basis functions $\left\{\phi_1(\vc x),\cdots, \phi_m(\vc x)\right\}$ and the measurements $\{u(\vc s_1),\cdots, u(\vc s_n)\}$, find the coefficients $c_i$ in~\eqref{eq:fin_series} so that $\tilde u \simeq u$.

\subsection{Finite-dimensional discretization}
\renewcommand{\arraystretch}{1.2}
\begin{table}[t]
	\centering
	\caption{The quantities and terminology used in this paper.}
	\begin{tabular}{|c |c |}
		\hline
		\textbf{Symbol} & \textbf{Description} \\ \hline\hline
		$N$ & High-fidelity Resolution \\ \hline
		$n$ & Number of Sensors \\ \hline
		$m$ & Number of Modes (Basis Functions) \\ \hline
		$\vc u\in\mathbb R^N$ & Unknown State\\ \hline
		$\vc y\in\mathbb R^n$ & Known Observations\\ \hline
		$S_n \in\mathbb R^{N\times n}$ & Selection Matrix, $\vc y=S_n^\top\vc u$\\ \hline
		$\Phi_m = [\pmb \phi_1|\cdots|\pmb \phi_m]\in\mathbb R^{N\times m}$ & Basis Matrix\\ \hline
		$\tilde{\vc u}\in\mathbb R^N$ & Reconstructed State\\ \hline
		$\hat{\vc u}\in\mathbb R^N$ & Orthogonal Reconstruction, $\hat{\vc u}=\Phi_m\Phi_m^\top \vc u$\\ \hline
		$\vc c\in\mathbb R^m$ & Expansion Coefficients, $\tilde{\vc u} = \Phi_m\vc c$\\ \hline
		$\mathbb D_{m,n} = \Phi_m (S_n^\top\Phi_m)^+S_n^\top\in \mathbb R^{N\times N}$ & DEIM Operator \\ \hline
	\end{tabular}
	\label{tab:glossary}
\end{table}

In practice, it is often sufficient to know the function $u$ over a dense enough grid $\mathcal G:=\{\vc x_i\}_{i=1}^N\subset \Omega$. We assume that $N$ is large enough so that $u(\vc x_i)$, $i=1,2,\cdots,N$ amounts to a high-resolution discretization of the function $u$.
This allows us to transform the problem from the infinite-dimensional function space to a finite-dimensional setting. We stack the function values over the high-resolution grid $\mathcal G$ into a vector 
\begin{equation}
\vc u = \begin{bmatrix}
u(\vc x_1),
u(\vc x_2),
\cdots,
u(\vc x_N)
\end{bmatrix}^\top\in\mathbb R^N.
\end{equation}
Note, however, that the function values are not known over the entire high-resolution grid $\mathcal G$. Instead of interpolating the function $u$ over the entire domain $\Omega$, we would like to find its values over the grid $\mathcal G$.

Similarly, we discretize the basis functions $\phi_i$ and define
\begin{equation}
\bphi_i = \begin{bmatrix}
\phi_i(\vc x_1),
\phi_i(\vc x_2),
\cdots,
\phi_i(\vc x_N)
\end{bmatrix}^\top\in\mathbb R^N, \quad i=1,2,\cdots, m.
\end{equation}
We also define the basis matrix $\Phi_m = [\bphi_1|\bphi_2|\cdots|\bphi_m]\in\mathbb R^{N\times m}$.
We assume that the basis vectors are orthonormalized, so that $\Phi_m^\top\Phi_m = \mathbb I_m$.
Finally, approximation~\eqref{eq:fin_series} is discretized to define
\begin{equation}
\tilde{\vc u} = \begin{bmatrix}
\tilde u(\vc x_1),
\tilde u(\vc x_2),
\cdots,
\tilde u(\vc x_N)
\end{bmatrix}^\top\in\mathbb R^N, 
\end{equation}
so that $\tilde{\vc u}=\Phi\vc c$, where $\vc c = [c_1,c_2,\cdots,c_m]^\top\in\mathbb R^m$ is the vector of unknown coefficients.

Recall that we seek to infer the coefficients $\vc c$ from the observational data $\{u(\vc s_1),\allowbreak u(\vc s_2),\cdots, u(\vc s_n)\}$ which are obtained at the sensor locations $\mathcal S = \{\vc s_k\}_{k=1}^n\allowbreak\subset \mathcal G$. We define the vector of observations $\vc y\in\mathbb R^n$ as
\begin{equation}
\vc y = \begin{bmatrix}
u(\vc s_1), &
u(\vc s_2),&
\cdots, &
u(\vc s_n)
\end{bmatrix}^\top + \pmb\varepsilon,
\end{equation}
where $\pmb\varepsilon\in\mathbb R^n$ is a vector of errors accounting for observational noise. In the following, we only consider the ideal case where the observational noise is negligible and set $\pmb\varepsilon=\vc 0$. Nonetheless, in section~\ref{sec:numerics} where our numerical results are presented, we study the effects of noise empirically.

It is convenient to envision the sensor measurements as a linear map from the full state $\vc u$ to the observations $\vc y$. To this end, we define the selection operator $S_n:\mathbb R^n\to\mathbb R^N$ which is an $N\times n$ matrix whose columns are a subset of the standard basis $\{\vc e_1,\vc e_2,\cdots,\vc e_N\}$ on $\mathbb R^N$. More specifically, let $\{i_1,i_2,\cdots,i_n\}$ denote a subset of the index set $\{1,2,\cdots,N\}$ such that $\vc x_{i_k}=\vc s_k$ for $k=1,2,\cdots, n$.
In other words, the distinct indices $i_k$ are chosen so that $\vc x_{i_k}$ coincides with the sensor locations $\vc s_k$. Consequently, defining the \emph{selection matrix} $S_n=[\vc e_{i_1}|\cdots|\vc e_{i_n}]$, we have
\begin{equation}
\vc y = S_n^\top \vc u.
\end{equation}
Therefore, $S_n^\top:\mathbb R^N\to\mathbb R^n$ returns the smaller vector $\vc y$ whose entries coincide with the entries of the full state vector $\vc u$ at the sensor locations $\vc x_{i_k}=\vc s_k$.
Conversely, $S_n\vc y$ returns a vector in $\mathbb R^N$ whose entries are zero except the ones corresponding to the sensor locations. In principle, $S_n\vc y$ is a reconstruction 
of the full state; however, this reconstruction is very inaccurate since it zeros out all the entries where no sensors are present.
The selection matrix $S_n$ plays an important role in DEIM. In particular, the optimal sensor placement can be equivalently expressed in terms of the selection matrix, as we describe in Section~\ref{sec:DEIM}. In the following, we often omit the subscripts $m$ and $n$, unless they are necessary for clarity. The quantities and terminology that are used repeatedly in this paper are summarized in Table~\ref{tab:glossary}.

\subsection{Vanilla DEIM}\label{sec:DEIM}
We refer to DEIM and Q-DEIM methods~\cite{Sorensen2010,Drmac2016} as vanilla DEIM to distinguish them from sparse DEIM introduced in section~\ref{sec:S-DEIM}.
Here, we first review vanilla DEIM for the special case where the number of sensors and modes are equal, $n=m$ (section~\ref{sec:equal}).
This is the setting in which vanilla DEIM is used most often.
Towards the end of this section, we discuss the implications of using vanilla DEIM with fewer sensors than modes, $n<m$ (section~\ref{sec:diff}).
DEIM was originally derived in a reduced-order modeling framework. Here, we present a different, but equivalent, derivation which is more natural for interpolation problems.

\subsubsection{Equal number of sensors and modes}\label{sec:equal}
One can seek the optimal coefficients $\vc c$ by solving the optimization problem, 
\begin{equation}
\min_{\vc c\in\mathbb R^m} \|\Phi\vc c-\vc u\|^2,
\label{eq:LS}
\end{equation}
where $\|\cdot\|$ denotes the Euclidean 2-norm. This problem has a unique solution $\vc c=\Phi^\top\vc u$, leading to the reconstruction $\hat{\vc u} = \Phi\Phi^\top \vc u$ which is the nearest point in $\R[\Phi]$ (range of $\Phi$) to the state $\vc u$. This motivates the following definition. 
\begin{defn}[Orthogonal reconstruction]
	For a basis matrix $\Phi\in\mathbb R^{N\times m}$ and a vector $\vc u\in\mathbb R^N$, we refer to the orthogonal projection $\hat{\vc u}=\Phi\Phi^\top \vc u$ as the orthogonal reconstruction of $\vc u$. We refer to the corresponding error $\mathcal E_m(\vc u)=\|\hat{\vc u} -\vc u\|$ as the truncation error.
\end{defn}

The orthogonal reconstruction $\hat{\vc u}$ is not computable from observations since it requires the knowledge of the entire state $\vc u$.
Note that the optimization problem~\eqref{eq:LS} does not use the observational data $\vc y$.
One can incorporate the observational data by adding them as constraints to the optimization problem, and solve
\begin{subequations}\label{eq:LS_const}
	\begin{gather}
	\min_{\vc c\in\mathbb R^m} \|\Phi\vc c-\vc u\|^2,\label{eq:LS_const_min}\\
	\mbox{such that} \quad S^\top \Phi \vc c = \vc y.
	\label{eq:LS_const_const}
	\end{gather}
\end{subequations}
Recalling that $\vc y = S^\top \vc u$, the solution to this constraint optimization problem is given by (see, e.g., Refs.~\cite{Bjorck1996,Lawson1995})
\begin{equation}
\vc c = \left[ \Phi^\top -\Phi^\top S (S^\top \Phi \Phi^\top S)^{-1}S^\top(\Phi\Phi^\top-\mathbb I_N)\right]\vc u.
\label{eq:lsp_const_sol}
\end{equation}
It is easy to verify that $S^\top \Phi\vc c = \vc y$ and therefore the reconstruction $\tilde{\vc u}=\Phi\vc c$ agrees with the observations $\vc y=S^\top \vc u$ at the sensor locations.

Unfortunately, expression~\eqref{eq:lsp_const_sol} is still not computable since it relies on the entire state $\vc u$.
However, in the special case where $S^\top\Phi$ is invertible, this expression simplifies to 
\begin{equation}
\vc c = (S^\top \Phi)^{-1}S^\top\vc u=(S^\top \Phi)^{-1}\vc y,
\label{eq:DEIM}
\end{equation}
which is now computable from known selection matrix $S$, basis matrix $\Phi$, and the observations $\vc y$. We note that when $S^\top\Phi$ is invertible, the constraint set~\eqref{eq:LS_const_const} contains the single point $\vc c =(S^\top \Phi)^{-1}\vc y$ and therefore the solution to the optimization problem is trivial. 

Expression~\eqref{eq:DEIM} leads to the vanilla DEIM reconstruction,
\begin{equation}\label{eq:vDEIM_rec}
\tilde{\vc u} = \Phi\vc c = \Phi (S^\top \Phi)^{-1}\vc y.
\end{equation}

The reconstruction error $\|\tilde{\vc u}-\vc u\|$ can be bounded from above with a multiple of the truncation error.

\begin{lem}[Vanilla DEIM error estimate~\cite{Sorensen2010}]\label{lem:vDEIM_err}
	Assume that the number of sensors and the number of modes are equal, $n=m$, and that the matrix $S^\top \Phi\in\mathbb R^{n\times n}$ is invertible. The vanilla DEIM reconstruction~\eqref{eq:vDEIM_rec} satisfies
	\begin{equation}\label{eq:vDEIM_err}
	\|\tilde{\vc u} - \vc u\| \leq \|(S^\top \Phi)^{-1}\|_2\ \mathcal E_{m}(\vc u),
	\end{equation}
	where $\mathcal E_{m}(\vc u)=\|\vc u-\hat{\vc u}\|$ is the truncation error, $\hat{\vc u} = \Phi\Phi^\top\vc u$ is the orthogonal reconstruction of $\vc u$, and $\|\cdot\|_2$ denotes the 2-norm (or spectral norm) of the matrix.
\end{lem}

The truncation error $\mathcal E_m$ decreases monotonically with the number of basis vectors $m$. 
Therefore, to decrease the vanilla DEIM reconstruction error, we ideally would need to use as many basis vector as possible. However, since $n=m$, the number of basis vectors $m$ is limited by the number of available sensors.

So far, we have said nothing about the choice of the selection matrix $S$ which determines the sensor locations. In vanilla DEIM, inequality~\eqref{eq:vDEIM_err} is used to inform the choice of the selection matrix. Ideally, the selection matrix should be chosen such that the prefactor 
$\|(S^\top \Phi)^{-1}\|_2$ is minimal and therefore the error upper bound~\eqref{eq:vDEIM_err} is as small as possible. Unfortunately, minimizing $\|(S^\top \Phi)^{-1}\|_2$ over all selection matrices is combinatorially hard and therefore intractable in most applications.
Therefore, one instead uses a greedy algorithm to obtain a suboptimal, yet acceptable, selection matrix $S$. The state-of-the-art uses QR factorization with column pivoting of $\Phi^\top$. Since sensor placement is not the focus of this paper, we refer to~\cite{Drmac2016} for further details and only state the QR algorithm using Matlab's syntax:
\begin{equation}\label{eq:QRalg}
\texttt{[\texttildelow ,\texttildelow,}\mathcal I\texttt{] = qr(}\Phi^\top\texttt{,`vector')}; \quad
\mathbb I\texttt{=eye(N)}; \quad
\texttt{S=}\mathbb I\texttt{(:,}\mathcal I(1:n)\texttt{)}.
\end{equation}
The set $\mathcal I(1:n)$ contains the indices $\{i_1,\cdots,i_n\}$ of the columns of the identity matrix $\mathbb I$ which form the selection matrix $S$.
The DEIM algorithm with QR factorization for sensor placement is often referred to as Q-DEIM. Here, we refer to both methods (DEIM and Q-DEIM) as vanilla DEIM to distinguish them from the S-DEIM method introduced in section~\ref{sec:S-DEIM}.

Now we turn our attention to some properties of the vanilla DEIM reconstruction~\eqref{eq:vDEIM_rec}. Recalling that $\vc y =S^\top \vc u$, 
we have $\tilde{\vc u} = \Phi (S^\top \Phi)^{-1} S^\top \vc u$. Defining the \emph{DEIM operator} $\mathbb D = \Phi (S^\top \Phi)^{-1} S^\top$, we have 
$\tilde {\vc u} = \mathbb D\vc u$. We note that the DEIM operator $\mathbb D$ is an oblique projection, parallel to the null space of $S^\top$, onto the range of $\Phi$~\cite{Drmac2018}. The following well-known result outlines some other properties of the DEIM operator.

\begin{lem}[Ref.~\cite{Drmac2018}]
	Consider the DEIM operator $\mathbb D = \Phi (S^\top \Phi)^{-1} S^\top$ with equal numbers of sensors and modes, $m=n$.
	We have
	\begin{enumerate}
		\item Interpolation property: $S^\top \mathbb D = S^\top$. 
		\item Projection property: $\mathbb D\Phi\Phi^\top = \Phi\Phi^\top$.
	\end{enumerate}
\end{lem}

The interpolation property implies that the observations $\tilde{\vc y} = S^\top \tilde{\vc u}$ corresponding to the DEIM reconstruction agree with the true observations $\vc y=S^\top \vc u$. In other words, the interpolation property implies that the reconstruction reproduces exact values at the sensor locations, although it can incur errors elsewhere.

The projection property implies that the vanilla DEIM reconstruction of the orthogonal projection is exact, $\mathbb D \hat{\vc u}=\hat{\vc u}$.
In particular, if $\vc u\in\R[\Phi]$, then its vanilla DEIM reconstruction returns the exact state since $\mathbb D\vc u=\vc u$.

\subsubsection{Different number of sensors and modes}\label{sec:diff}
In this section, we comment on the application of vanilla DEIM when the number of sensors and modes are not equal, $n\neq m$.
In this case, the matrix $S_n^\top \Phi_m$ is no longer square and therefore its inverse is not well-defined. Instead, one replaces the inverse with the Moore--Penrose pseudo inverse $(S_n^\top \Phi_m)^+$, so that the expansion coefficients are given by $\vc c =(S_n^\top\Phi_m)^+\vc y$. Note that this coefficient vector is the \emph{minimum-norm} solution to Eq.~\eqref{eq:LS_const_const}.
The corresponding vanilla DEIM reconstruction is given by $\tilde{\vc u} = \mathbb D_{m,n} \vc u$, where the DEIM operator is now defined by $\mathbb D_{m,n} = \Phi_m(S_n^\top\Phi_m)^+S_n^\top$. To obtain the selection matrix $S_n$, one uses the same QR factorization with column pivoting as outlined in~\eqref{eq:QRalg}.

In practice, the matrix $S_n^\top\Phi_m$ often has full rank. In keeping with previous studies~\cite{Sorensen2010,Drmac2018}, we make this assumption throughout this paper.
\begin{ass}\label{ass:full-rank}
	For a given basis matrix $\Phi_m$, let the selection matrix $S_n$ be determined by the Q-DEIM algorithm~\eqref{eq:QRalg}. We assume that $S_n^\top\Phi_m$ has full rank, i.e., 
	$\mbox{rank}[S_n^\top\Phi_m]=\min\{n,m\}$.
\end{ass}

\begin{rem}[See Ref.~\cite{Drmac2018}]\label{rem:properties}
	Unfortunately, when $n\neq m$, some of the nice properties of vanilla DEIM fall apart:
	\begin{enumerate}
		\item If $n<m$ (underdetermined, fewer sensor than modes), the interpolation property, $S_n^\top \mathbb D_{m,n}=S_n^\top$, still holds. However, the projection property fails, i.e., $\mathbb D_{m,n}\Phi_m\Phi_m^\top \neq \Phi_m\Phi_m^\top$. Furthermore, the error upper bound~\eqref{eq:vDEIM_err} is no longer valid.
		\item If $n>m$ (overdetermined, more sensor than modes), the interpolation property fails, i.e., $S_n^\top \mathbb D_{m,n}\neq S_n^\top$. But the projection property and the error upper bound are still valid.
	\end{enumerate}
	These statements are summarized in Table~\ref{tab:prop}.
\end{rem}

\begin{table}[h!]
	\centering
	\caption{Vanilla DEIM properties with various number of sensors $n$ and modes $m$.}
	\begin{tabular}{|c |c |c | c |}
		\hline
		&  $n=m$ &  $n>m$ &  $n<m$ \\ \hline
		Interpolation property & $\checkmark$ & $\times$ & $\checkmark$ \\ \hline
		Projection property & $\checkmark$ & $\checkmark$ & $\times$ \\ \hline
		Error upper bound & $\checkmark$ & $\checkmark$ & $\times$ \\ \hline
	\end{tabular}
	\label{tab:prop}
\end{table}

Our main focus in this paper is on the underdetermined case where the number of sensors is smaller than the number of modes, $n<m$. In this regime, under Assumption~\ref{ass:full-rank}, the pseudo inverse $(S_n^\top\Phi_m)^+$ is in fact a right inverse: $(S_n^\top\Phi_m)(S_n^\top\Phi_m)^+=\mathbb I_n$. In section~\ref{sec:S-DEIM}, where we introduce sparse DEIM, we show that the appropriate error upper bound is more complicated than Eq.~\eqref{eq:vDEIM_err} and involves the null space of $S_n^\top \Phi_m$.

\begin{figure}
	\centering
	\includegraphics[width=.9\textwidth]{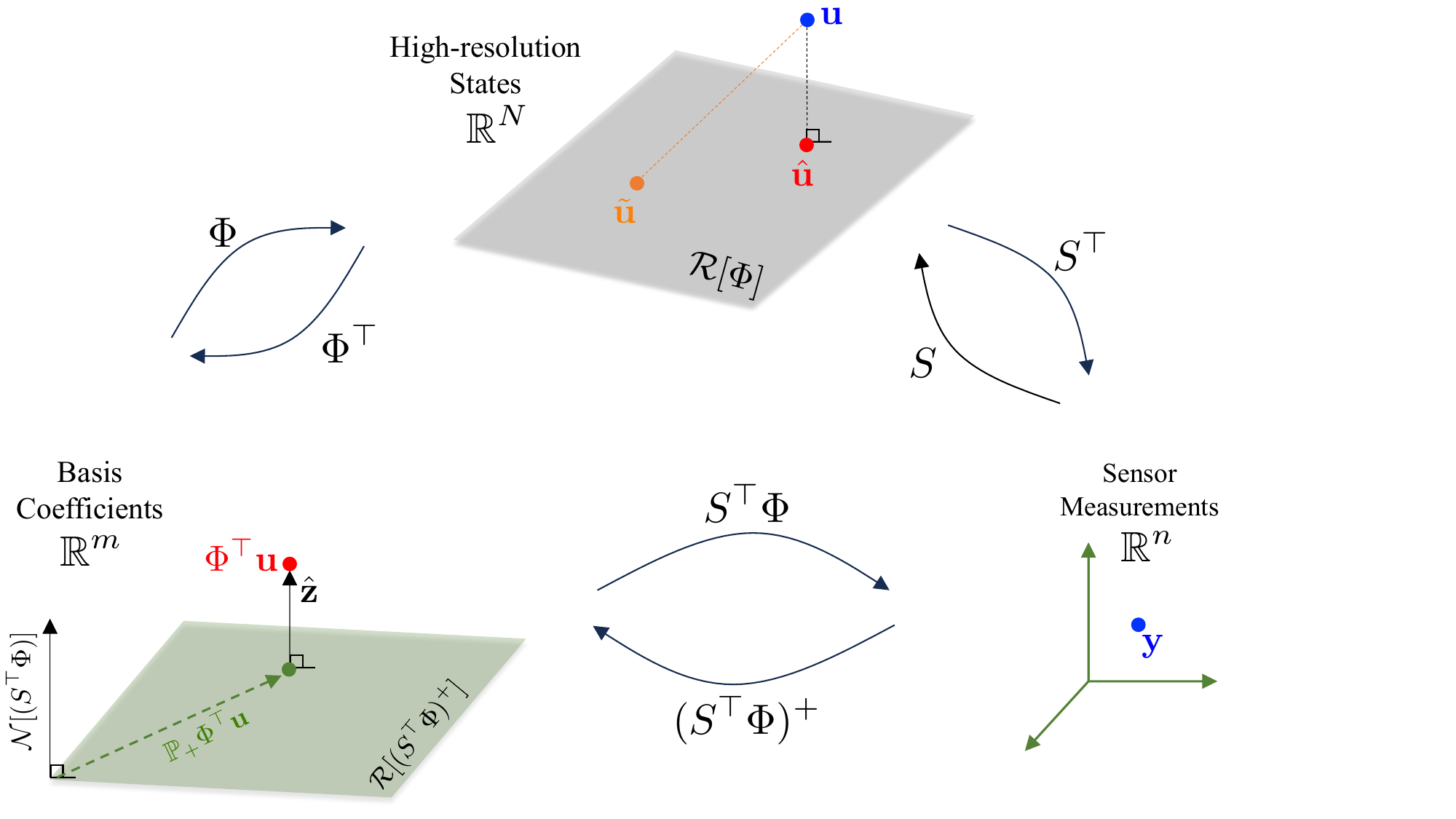}
	\caption{Mappings between the high-resolution space $\mathbb R^N$, the measurement space $\mathbb R^n$, and the coefficient space $\mathbb R^m$. The matrix $\mathbb P_+$ is shorthand for $\mathbb P_{\R[(S^\top\Phi)^+]}$, i.e., orthogonal projection onto the range of $(S^\top\Phi)^+$. Note that the mappings in this figure do not necessarily commute.}
\end{figure}

\section{Sparse DEIM}\label{sec:S-DEIM}
Vanilla DEIM uses the expansion coefficients $\vc c=(S^\top \Phi)^+\vc y$ which are the \emph{minimum-norm} solution to the least squares problem,
\begin{equation}\label{eq:LSP}
\min_{\vc c\in\mathbb R^m}\|S^\top \Phi\vc c - \vc y\|^2. 
\end{equation}
However, the general solution to this least squares problem is given by
\begin{equation}
\vc c(\vc z) = (S^\top \Phi)^+\vc y +\vc z,\quad \forall \vc z\in \Null [S^\top \Phi],
\end{equation}
where $\Null$ denotes the null space.
Note that $\vc c(\vc 0)$ corresponds to the vanilla DEIM reconstruction. If the number of sensors are equal to or greater than the number of modes, $n\geq m$, then we have $\Null[S^\top \Phi]=\{\vc 0\}$ and the vanilla DEIM reconstruction is the only choice. However, when the number of sensors is smaller than the number of modes, $n<m$, the null space $\Null[S^\top\Phi]$ is an $(m-n)$-dimensional subspace of $\mathbb R^m$ and $\vc z\neq \vc 0$ is a possibility. This motivates the following definition. 
\begin{defn}[Sparse DEIM reconstruction]\label{def:S-DEIM}
	Assume there are fewer sensors than modes, $n<m$. For any nonzero $\vc z\in \Null[S^\top\Phi]$, we refer to 
	\begin{equation}\label{eq:S-DEIM}
	\tilde{\vc u}(\vc z) = \Phi\vc c(\vc z) = \Phi(S^\top \Phi)^+\vc y +\Phi\vc z=\mathbb D \vc u + \Phi\vc z,
	\end{equation}
	as the \emph{sparse DEIM reconstruction} of $\vc u$, or S-DEIM for short. We refer to $\vc z\in \Null[S^\top\Phi]$ as a \emph{kernel vector}.
\end{defn}

The main question is whether there exists $\vc 0\neq\vc z\in\Null[S^\top \Phi]$ such that the corresponding S-DEIM reconstruction $\tilde{\vc u}(\vc z)$ is a better approximation of $\vc u$ than the vanilla DEIM reconstruction $\tilde{\vc u}(\vc 0)$?

We first note that, under Assumption~\ref{ass:full-rank}, the addition of $\vc z\neq\vc 0$ does not change the value of the reconstructed field at the sensor locations. More precisely, 
\begin{equation}
S^\top \tilde{\vc u}(\vc z) = S^\top \Phi(S^\top \Phi)^+ \vc y+S^\top\Phi\vc z= \vc y,\quad \forall \vc z\in\Null[S^\top \Phi],
\end{equation}
which implies that any S-DEIM reconstruction agrees with the data from sensor measurements, i.e., S-DEIM respects the interpolation property. As a result, the available sensor data cannot be immediately used to inform the best choice of the kernel vector $\vc z$. Here, we defer the choice of $\vc z$ to sections~\ref{sec:two-stage} and~\ref{sec:DAS-DEIM}, and first derive error estimates for the S-DEIM reconstruction with an arbitrary kernel vector $\vc z\in\Null[S^\top \Phi]$.

Before presenting our main error estimate in Theorem~\ref{thm:error}, we need some preliminary results. 
\begin{lem}\label{lem:orth_decomp}
	Range of $(S^\top \Phi)^+$ is the orthogonal complement of the null space of $S^\top \Phi$, i.e.,
	\begin{equation}
	\mathbb R^m = \R [(S^\top \Phi)^+] \oplus \Null[S^\top \Phi] \quad \mbox{and}\quad
	\R [(S^\top \Phi)^+]\perp \Null[S^\top \Phi]
	\end{equation}
\end{lem}
\begin{proof}
	This is a well-known consequence of the fundamental theorem of linear algebra (see, e.g., Ref.~\cite{Strang1993}) and the fact that $\R[(S^\top \Phi)^+]=\R[(S^\top \Phi)^\top]$.
\end{proof}

\begin{figure}
	\centering
	\includegraphics[width=.5\textwidth]{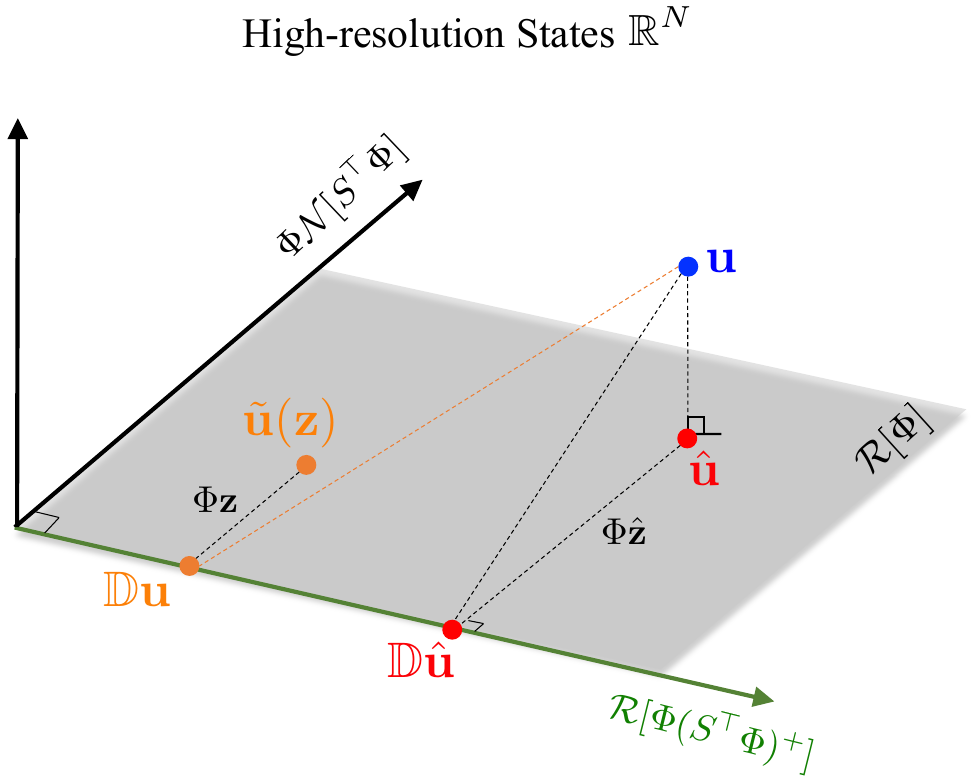}
	\caption{Orthogonal decomposition of $\R[\Phi]$ with fewer sensors than modes, $n<m$. Dashed black lines represent orthogonal projections, whereas the dashed orange line represents an oblique projection.}
	\label{fig:RangePhi}
\end{figure}

As depicted in figure~\ref{fig:RangePhi}, Lemma~\ref{lem:orth_decomp} implies that the range of $\Phi$ can be decomposed into two orthogonal subspaces as follows. 
\begin{cor}\label{cor:orth_decomp}
	The range of $\Phi$ admits the following orthogonal decomposition,
	\begin{equation}
	\R[\Phi] = \R [\Phi(S^\top \Phi)^+] \oplus \Phi\Null[S^\top \Phi] \quad \mbox{and}\quad
	\R [\Phi(S^\top \Phi)^+]\perp \Phi\Null[S^\top \Phi],
	\end{equation}
	where $\Phi\Null[S^\top \Phi]=\{\vc u\in\mathbb R^N: \vc u = \Phi\vc z,\ \vc z\in\Null[S^\top \Phi]\}$.
\end{cor}
\begin{proof}
	This follows from the fact that the columns of $\Phi$ are orthonormal. In particular, it is injective and it preserves the Euclidean inner product ($\Phi^\top \Phi=\mathbb I_m$).
\end{proof}

\begin{defn}[Kernel matrix]
	Let $\{\vc z_1,\vc z_2,\cdots,\vc z_{m-n}\}\subset \mathbb R^m$ be an orthonormal basis for $\Null[S^\top \Phi]$. We refer to the matrix $Z=[\vc z_1|\vc z_2|\cdots|\vc z_{m-n}]$ as a \emph{kernel matrix}.
\end{defn}

The following lemma gives the best choice of the kernel vector $\vc z\in\Null[S^\top \Phi]$ in terms of the kernel matrix $Z$.
\begin{lem}\label{lem:best_z}
	The optimization problem,
	\begin{equation}\label{eq:opt_best_z}
	\hat{\vc z}=\argmin_{\vc z \in \Null[S^\top \Phi]} \|\tilde{\vc u}(\vc z) - \vc u\|^2,
	\end{equation}
	has a unique solution given by $\hat{\vc z} = ZZ^\top \Phi^\top \vc u,$
	where $Z$ is a kernel matrix.
\end{lem}
\begin{proof}
	Consider a kernel matrix $Z$ and note that, for any $\vc z\in \Null[S^\top \Phi]$, there exists a unique $\bxi\in\mathbb R^{m-n}$ such that $\vc z = Z\bxi$. Therefore, optimization problem~\eqref{eq:opt_best_z} can be equivalently written as
	\begin{equation}
	\min_{\bxi \in \mathbb R^{m-n}} \|\Phi(S^\top \Phi)^+\vc y +\Phi Z\bxi- \vc u\|^2.
	\end{equation} 
	Since $\Phi Z$ has full column rank, the solution is unique. Furthermore, since $(\Phi Z)^+ = Z^\top \Phi^\top$, the solution is given by
	\begin{equation}
	\hat\bxi = Z^\top\Phi^\top \vc u - Z^\top \Phi^\top \Phi(S^\top \Phi)^+\vc y.
	\end{equation}
	Next we use the fact that $\Phi^\top \Phi=\mathbb I_m$ and $Z^\top(S^\top \Phi)^+=0$. The latter identity follows from Lemma~\ref{lem:orth_decomp}:   $\R(Z) = \Null[S^\top \Phi]\perp \R[(S^\top \Phi)^+]$. Therefore, we have $\hat\bxi = Z^\top\Phi^\top \vc u$ and $\hat{\vc z} = ZZ^\top\Phi^\top \vc u$. Note that, although the kernel matrix $Z$ is not necessarily unique, $ZZ^\top$ is the orthogonal projection onto $\Null[S^\top \Phi]$ and therefore $\hat{\vc z}$ is unique.
\end{proof}

Lemma~\ref{lem:best_z} gives the best choice of the kernel vector $\vc z$ which minimizes the error of the S-DEIM reconstruction $\tilde{\vc u}(\vc z)$. Unfortunately, the optimal kernel vector $\hat{\vc z} = ZZ^\top \Phi^\top \vc u$ is not computable from observations since it requires the entire state $\vc u$. In section~\ref{sec:DAS-DEIM}, we propose a computable method for estimating the optimal kernel vector. But first, we state our main error estimates of the S-DEIM reconstruction for arbitrary kernel vectors.

\begin{thm}[Sparse DEIM error]\label{thm:error}
	Let the number of sensors be smaller than the number of modes, $n<m$.
	For $\vc u\in\mathbb R^N$, let $\hat{\vc u} = \Phi\Phi^\top \vc u$ denote its orthogonal projection onto $\R[\Phi]$. 
	Let $\tilde{\vc u}(\vc z)$ denote the S-DEIM reconstruction for an arbitrary kernel vector $\vc z\in \Null[S^\top \Phi]$, and $\hat{\vc z} = ZZ^\top \Phi^\top \vc u$ as in Lemma~\ref{lem:best_z}. Then
	\begin{equation}
	\|\vc u - \tilde{\vc u}(\vc z)\|^2 = \|\vc u - \hat{\vc u}\|^2+\|\mathbb D(\vc u - \hat{\vc u})\|^2+\|\hat{\vc z} - \vc z\|^2,\quad \forall\vc z\in \Null[S^\top \Phi].
	\end{equation}
\end{thm}
\begin{proof}
	First we write $\vc u-\tilde{\vc u}(\vc z)=\vc u-\hat{\vc u}+\hat{\vc u}-\tilde{\vc u}(\vc z)$ where $\hat{\vc u}=\Phi\Phi^\top \vc u$ and $\tilde{\vc u}(\vc z) = \Phi\vc c(\vc z)$ both belong to $\R[\Phi]$; therefore, $\hat{\vc u}-\tilde{\vc u}(\vc z)\in\R[\Phi]$. On the other hand, $\vc u-\hat{\vc u}= (\mathbb I-\Phi\Phi^\top)\vc u$ is orthogonal to $\R[\Phi]$. Therefore, 
	\begin{equation}
	\|\vc u-\tilde{\vc u}(\vc z)\|^2 = \|\vc u-\hat{\vc u}+\hat{\vc u}-\tilde{\vc u}(\vc z)\|^2  = \| \vc u-\hat{\vc u}\|^2+\|\hat{\vc u}-\tilde{\vc u}(\vc z)\|^2.
	\end{equation}
	
	Next we focus on $\hat{\vc u}-\tilde{\vc u}(\vc z)$. Recall from Lemma~\ref{lem:orth_decomp} that $\mathbb R^m = \Null [S^\top\Phi]\oplus \R[(S^\top\Phi)^+]$ and that the two subspaces are orthogonal.
	Write $\Phi^\top \vc u\in\mathbb R^m$ as its unique orthogonal decomposition, 
	\begin{equation}
	\Phi^\top \vc u = \mathbb P_{\R[(S^\top\Phi)^+]}\Phi^\top\vc u + \hat{\vc z},
	\label{eq:PhiTu}
	\end{equation}
	where $\hat{\vc z}\in \Null [S^\top\Phi]$ is unique. Note that the orthogonal projection $ \mathbb P_{\R[(S^\top\Phi)^+]}$ onto $\R[(S^\top\Phi)^+]$ can be written as
	\begin{equation}
	\mathbb P_{\R[(S^\top\Phi)^+]} = (S^\top\Phi)^+(S^\top\Phi),
	\label{eq:PSTPhipinv}
	\end{equation}
	which implies
	\begin{equation}
	\hat{\vc u}=\Phi \Phi^\top \vc u = \Phi (S^\top\Phi)^+(S^\top\Phi)\Phi^\top \vc u + \Phi \hat{\vc z} = \mathbb D \Phi \Phi^\top \vc u + \Phi \hat{\vc z} =  \mathbb D \hat{\vc u} + \Phi \hat{\vc z}.
	\end{equation}
	As a side note, we mention that, if $n=m$ and $S^\top\Phi$ is invertible, then $\mathbb D\Phi\Phi^\top = \Phi\Phi^\top$. However, in our case where $n<m$, the identity does not hold (cf. Remark~\ref{rem:properties}).
	
	Returning to the proof and recalling that $\tilde{\vc u}(\vc z)=\Phi\vc c(\vc z) = \mathbb D\vc u + \Phi\vc z$, we have
	\begin{align}
	\|\hat{\vc u} - \tilde{\vc u}(\vc z)\|^2 & = \|\mathbb D \hat{\vc u}  - \mathbb D\vc u + \Phi\hat{\vc z} - \Phi\vc z\|^2\nonumber\\
	& =  \|\mathbb D \hat{\vc u} - \mathbb D\vc u\|^2 + \|\Phi \hat{\vc z}-\Phi\vc z\|^2\nonumber\\
	& =  \|\mathbb D \hat{\vc u} - \mathbb D\vc u\|^2 + \|\hat{\vc z}-\vc z\|^2,
	\end{align}
	where we used the fact that $\Phi\hat{\vc z} - \Phi\vc z\in\Phi\Null[S^\top\Phi]$ and $\mathbb D \Phi\Phi^\top \vc u - \mathbb D\vc u\in\R[\Phi(S^\top\Phi)^+]$ are orthogonal (cf. Corollary~\ref{cor:orth_decomp}).
	
	It remains to show that $\hat{\vc z} = ZZ^\top\Phi^\top\vc u$. Recall from equations~\eqref{eq:PhiTu} and~\eqref{eq:PSTPhipinv} that
	\begin{equation}
	\Phi^\top \vc u = (S^\top\Phi)^+S^\top \Phi\Phi^\top\vc u + \hat{\vc z}.
	\end{equation}
	For any $\hat{\vc z}\in\Null[S^\top \Phi]$ there exists a unique $\bxi\in\mathbb R^{n-m}$ such that $\hat{\vc z} = Z\bxi$. Substituting this in the above equation, and noting that $\Null[S^\top \Phi]$ is orthogonal to $\R[(S^\top\Phi)^+]$, we obtain $\bxi = Z^\top \Phi^\top \vc u$.
\end{proof}

Theorem~\ref{thm:error} has a number of important consequences. For instance, if the true state $\vc u$ belongs to the range of $\Phi$, then there exists an exact S-DEIM reconstruction.
\begin{cor}\label{cor:exactReconst}
	Let the number of sensors be smaller than the number of modes, $n<m$. If $\vc u\in \R(\Phi)$, then the S-DEIM reconstruction with $\vc z = \hat{\vc z} = ZZ^\top\Phi^\top\vc u$ is exact:
	$\tilde{\vc u}(\hat{\vc z}) = \vc u$.
\end{cor}
\begin{proof}
	Choose $\vc z = \hat{\vc z}$ in Theorem~\ref{thm:error}. Since $\vc u\in\R[\Phi]$, we have $\hat{\vc u} = \vc u$ and therefore $\|\vc u - \tilde{\vc u}(\hat{\vc z})\|=0$.
\end{proof}

\begin{rem}
	Corollary~\ref{cor:exactReconst} shows that, if $\vc u\in\R[\Phi]$, then the S-DEIM reconstruction is exact. Note, however, that when $n<m$ the vanilla DEIM reconstruction $\mathbb D\vc u$ is not exact, $\mathbb D\vc u\neq \vc u$, even if $\vc u\in\R[\Phi]$. This is due to the fact that the vanilla DEIM reconstruction fails to preserve the projection property as discussed in Remark~\ref{rem:properties}, i.e., $\mathbb D\Phi\Phi^\top\neq \Phi\Phi^\top$.
	
	Although $\mathbb D$ is not an orthogonal projection, it is straightforward to verify that $\mathbb D\Phi\Phi^\top = \Phi(S^\top\Phi)^+(S^\top\Phi)\Phi^\top$ is an orthogonal projection onto $\R[\Phi(S^\top\Phi)^+]$. As a result, even if $\hat{\vc u}\neq \vc u$, $\mathbb D\hat{\vc u} = \mathbb D\Phi\Phi^\top\vc u$ is the orthogonal projection of $\vc u$ onto $\R[\Phi(S^\top\Phi)^+]$. Figure~\ref{fig:RangePhi} shows this schematically.
\end{rem}

Theorem~\ref{thm:error} also leads to an upper bound for the S-DEIM reconstruction error. This upper bound is a generalization of Eq.~\eqref{eq:vDEIM_err} to the case where the number of sensors is smaller than the number of modes.
\begin{cor}[Error estimate for S-DEIM]\label{cor:SDEIM_ub}
	Let the number of sensors be smaller than the number of modes, $n<m$. The S-DEIM reconstruction satisfies,
	\begin{equation}\label{eq:SDEIM_ub}
	\|\tilde{\vc u}(\vc z)-\vc u\|\leq \|(S^\top\Phi)^+\|_2\, \mathcal E_m(\vc u)+\|\vc z-\hat{\vc z}\|,\quad \forall\vc z\in\Null[S^\top\Phi],
	\end{equation}
	where $\mathcal E_m(\vc u)= \|\hat{\vc u}-\vc u\|$ is the truncation error and $\hat{\vc z} = ZZ^\top\Phi^\top\vc u$ is the optimal kernel vector.
\end{cor}
\begin{proof}
	We recall from the proof of Theorem~\ref{thm:error} that $\tilde{\vc u}(\vc z) = \mathbb D\vc u+\Phi\vc z$ and $\hat{\vc u} = \mathbb D \hat{\vc u} + \Phi \hat{\vc z}$. Therefore, we have $\tilde{\vc u}(\vc z) - \vc u = (\mathbb I-\mathbb D)(\hat{\vc u} - \vc u) + \Phi(\vc z- \hat{\vc z})$. Using the triangle inequality, we obtain
	\begin{align}
	\|\tilde{\vc u}(\vc z) - \vc u\| & \leq \| (\mathbb I-\mathbb D)(\hat{\vc u} - \vc u)\|+\| \Phi(\vc z- \hat{\vc z}) \|\nonumber\\
	& \leq \|\mathbb D\|_2\, \mathcal E_m(\vc u) + \| \vc z- \hat{\vc z}\|,
	\end{align}
	where in the last inequality, we used the fact that $\mathbb D$ is a non-zero and non-identity projection, $\mathbb D^2=\mathbb D$. Any such projection satisfies $\|\mathbb D\|_2 = \|\mathbb I -\mathbb D\|_2$~\cite{Ipsen1995,Szyld2006}. Finally, we use the fact that the DEIM operator 
	$\mathbb D = \Phi(S^\top\Phi)^+S^\top$ satisfies $\|\mathbb D\|_2 \leq \|(S^\top\Phi)^+\|_2$ since $\|\Phi\|_2=\|S^\top\|_2=1$.
\end{proof}

When the number of sensors is equal to the number of modes, $n=m$, we have $\Null[S^\top\Phi]=\{\vc 0\}$ and therefore $\vc z=\hat{\vc z} =\vc 0$.
As a result, the last term on the right-hand side of Eq.~\eqref{eq:SDEIM_ub} vanishes. In this special case, the error upper bound coincides with the one in Lemma~\ref{lem:vDEIM_err} for vanilla DEIM. However, when $n<m$, the correct upper bound is the one given in Eq.~\eqref{eq:SDEIM_ub} involving the difference between the kernel vector $\vc z$ and its optimal value $\hat{\vc z}$.

\begin{rem}\label{rem:prefactor}
	The error upper bound~\eqref{eq:SDEIM_ub} highlights an important advantage of S-DEIM over underdetermined vanilla DEIM ($n<m$). Assume that a kernel vector $\vc z$ is chosen such that the kernel error is small, $\|\vc z-\hat{\vc z}\|<\epsilon$. By increasing the number of modes $m$, the S-DEIM reconstruction error will eventually become smaller than $\epsilon$. This is because the truncation error $\mathcal E_m(\vc u)$ is a monotonically decreasing function of $m$. Furthermore, for a fixed number of sensors, the spectral norm $\|(S_n^\top\Phi_m)^+\|_2$ decreases as $m$ increases (This follows from interlacing inequalities for singular values; see, e.g., Theorem 7.3.9 of~\cite{Horn1985}). Such a small error cannot be guaranteed for underdetermined vanilla DEIM since $\vc z=\vc 0$ and therefore the upper bound is always larger than $\|\hat{\vc z}\|$.
\end{rem}


In this section, we largely sidestepped the choice of the kernel vector $\vc z$. In section~\ref{sec:DAS-DEIM}, we propose a computable method for obtaining an optimal kernel vector $\vc z$ which minimizes the S-DEIM reconstruction error. But first, in section~\ref{sec:two-stage} below, we discuss a promising idea for selecting the kernel vector and show that it unfortunately fails to lead to any improvement over vanilla DEIM.

\subsection{Withholding measurements would not work}\label{sec:two-stage}
Assume that $n$ distinct sensor measurements are available. We split these measurements into two batches with $n_1$ and $n_2$ measurement in each so that $n_1+n_2=n$.
We withhold the $n_2$ measurements at first and use $n_1$ sensor measurements to obtain the S-DEIM reconstruction $\tilde{\vc u}(\vc z)$ with an unknown kernel vector $\vc z$. Then the remaining $n_2$ withheld measurements are used to compute the kernel vector $\vc z$.

To describe the method more precisely, consider a selection matrix $S\in\mathbb R^{N\times n}$ which identifies the position of $n$ distinct sensors. This selection matrix, for instance, can be obtained by the column pivoted QR factorization~\eqref{eq:QRalg}. We split the selection matrix into two blocks, $S = [S_1\quad S_2]$, where $S_1\in\mathbb R^{N\times n_1}$ corresponds to the first batch of $n_1$ sensors and $S_2\in\mathbb R^{N\times n_2}$ corresponds to the second batch. Note that we use the shorthand notation $S_i$ for $S_{n_i}$. Furthermore, we assume, without loss of generality, that the first batch corresponds to the first $n_1$ columns of $S$. We write the measurements $\vc y$ as
\begin{equation}
\vc y = \begin{bmatrix}
\vc y_1 \\ \vc y_2
\end{bmatrix},
\end{equation}
where $\vc y_1\in\mathbb R^{n_1}$ denotes the measurements from the first batch of sensors and $\vc y_2\in\mathbb R^{n_2}$ denotes the remaining measurements. Note that $\vc y_i=S_i^\top \vc u$, $i=1,2$, where $\vc u$ is the unknown true state.

With this set-up, and only using the first $n_1$ sensors, the S-DEIM coefficients are given by 
\begin{equation}
\vc c(\vc z) = (S_1^\top \Phi)^+\vc y_1 + \vc z,
\end{equation}
where $\vc z\in \Null [S_1^\top\Phi]$ is the kernel vector to be determined. The reconstructed state is given by $\tilde {\vc u}(\vc z) = \Phi\vc c(\vc z)$. 
In order to find the kernel vector $\vc z$, we use the second batch of sensors. More precisely, we determine $\vc z$ by requiring $S_2^\top\tilde{\vc u}(\vc z)$, i.e., the reconstructed state measured at the second batch of sensors, to be as close to the true measurements $\vc y_2$ as possible.
In other words, we solve the optimization problem,
\begin{equation}
\mmin_{\vc z\in \Null[S_1^\top\Phi]} \|S_2^\top \Phi \vc c(\vc z) - \vc y_2\|^2,
\label{eq:min_y2}
\end{equation}
where the overline indicates the minimum-norm solution.
As the following theorem states, the solution to~\eqref{eq:min_y2} coincides with the solution to vanilla DEIM reconstruction~\eqref{eq:vDEIM_rec} using all sensors at once.
\begin{thm}[Equivalence of vanilla DEIM and two-stage S-DEIM]
	Consider a S-DEIM reconstruction using $n_1$ sensors, $\tilde{\vc u}(\vc z) = \Phi\vc c(\vc z)$ where $\vc c(\vc z)=(S_1^\top \Phi)^+\vc y_1 + \vc z$. Determine the kernel vector $\vc z\in\Null[S_1^\top \Phi]$ as the solution of~\eqref{eq:min_y2}.
	Under Assumption~\ref{ass:full-rank}, $\tilde{\vc u}(\vc z)$ coincides with the vanilla DEIM reconstruction $\Phi(S^\top\Phi)^+\vc y$.
\end{thm}
\begin{proof}
	We first write the minimization problem~\eqref{eq:min_y2} more explicitly. Let $Z_1\in\mathbb R^{m\times (m-n_1)}$ denote a matrix whose columns form an orthonormal basis for $\Null [S_1^\top \Phi]$. Any $\vc z\in\Null [S_1^\top \Phi]$ can be written uniquely as $\vc z = Z_1\bxi$ for some $\bxi\in\mathbb R^{m-n_1}$.
	Therefore, optimization problem~\eqref{eq:min_y2} can be written equivalently as
	\begin{equation}
	\mmin_{\bxi\in\mathbb R^{m-n_1}}\|S_2^\top \Phi(S_1^\top \Phi)^+\vc y_1 + S_2^\top \Phi Z_1\bxi-\vc y_2\|^2
	\label{eq:min_y2_v2}
	\end{equation}
	
	Now we turn our attention to vanilla DEIM. Recall that vanilla DEIM corresponds to the minimum-norm minimizer of the cost function
	\begin{equation}
	J(\vc c) = \|S^\top \Phi \vc c- \vc y\|^2.
	\label{eq:J_vDEIM}
	\end{equation}
	First, we write this cost function in a slightly different but equivalent form. Note that, for all $\vc c\in\mathbb R^m$, we have
	\begin{equation}
	S^\top \Phi\vc c = \begin{bmatrix}
	S_1^\top \Phi\vc c \\
	S_2^\top \Phi\vc c
	\end{bmatrix}, 
	\end{equation}
	and therefore, 
	\begin{equation}
	J(\vc c) = \|S_1^\top \Phi \vc c - \vc y_1\|^2+\|S_2^\top \Phi \vc c - \vc y_2\|^2.
	\label{eq:split_cost}
	\end{equation}
	Next we decompose $\vc c$ into its orthogonal projections onto the subspaces $\R[(S_1^\top \Phi)^+]$ and $\Null [S_1^\top \Phi]$. Recall that these are subspaces of $\mathbb R^m$ and are orthogonal complements of each other. Therefore, for any $\vc c\in\mathbb R^m$, there exist $\etab\in\mathbb R^{n_1}$ and $\bxi\in\mathbb R^{m-n_1}$ such that
	$\vc c = (S_1^\top \Phi)^+\etab + Z_1\bxi$.
	Recalling that 	$S_1^\top\Phi(S_1^\top \Phi)^+ = \mathbb I_{n_1}$, we have
	\begin{subequations}
		\begin{equation}
		S_1^\top \Phi \vc c = \etab,
		\end{equation}
		\begin{equation}
		S_2^\top \Phi \vc c = S_2^\top \Phi(S_1^\top \Phi)^+\etab + S_2^\top \Phi Z_1\bxi.
		\end{equation}
	\end{subequations}
	Substituting these in~\eqref{eq:split_cost}, we obtain
	\begin{equation}
	J(\vc c) = J(\etab,\bxi) = \|\etab - \vc y_1\|^2 + \|S_2^\top \Phi(S_1^\top \Phi)^+\etab + S_2^\top \Phi Z_1\bxi-\vc y_2\|^2.
	\label{eq:J_vDEIM2}
	\end{equation}
	Therefore, the vanilla DEIM optimization $\mmin_{\vc c\in\mathbb R^m} J(\vc c)$ can be written equivalently as
	\begin{equation}
	\mmin_{\etab\in\mathbb R^{n_1}, \bxi\in\mathbb R^{m-n_1}} J(\etab,\bxi).
	\label{eq:opt_vDEIMxi}
	\end{equation} 
	We denote the minimizer by $(\etab_\ast,\bxi_\ast)$. 
	Now let $\vc c_\ast= (S^\top\Phi)^+\vc y$ denote the vanilla DEIM minimizer of $J(\vc c)$.
	Recall that the vanilla DEIM cost function attains $J(\vc c_\ast) = 0$ at its minimum. Since $J(\vc c) = J(\etab,\bxi)$, we must also have $J(\etab_\ast,\bxi_\ast)=0$. Examining equation~\eqref{eq:J_vDEIM2}, for $J(\etab_\ast,\bxi_\ast) =0$ to be achieved, it is necessary that $\etab_\ast=\vc y_1$. Therefore, optimization problem~\eqref{eq:opt_vDEIMxi} is equivalent to 
	\begin{equation}
	\mmin_{\bxi\in\mathbb R^{m-n_1}} \|S_2^\top \Phi(S_1^\top \Phi)^+\vc y_1 + S_2^\top \Phi Z_1\bxi-\vc y_2\|^2,
	\end{equation}
	where we substituted $\etab = \vc y_1$.
	This optimization problem is equivalent to vanilla DEIM and identical to the optimization problem~\eqref{eq:min_y2_v2}. Therefore, the two-stage S-DEIM is equivalent to the vanilla DEIM reconstruction.
\end{proof}

\section{Data assimilation with S-DEIM}\label{sec:DAS-DEIM}
So far we have not specified how the states $\vc u$ are generated. In this section, we assume that the states are generated by a continuous-time dynamical system. In this case, we propose a method for estimating the unknown kernel vector of S-DEIM using the observational time series $\vc y(t)$. The resulting method, which we refer to as the \emph{Data Assimilated Sparse DEIM} (DAS-DEIM), amounts to a fast data assimilation method for dynamical systems.

Consider a system of ODEs, possibly arising from the spatial discretization of a partial differential equation,
\begin{equation}
\frac{\id\vc u}{\id t} = \vc f(\vc u),\quad \vc u(0)=\vc u_0,
\label{eq:ODE}
\end{equation}
where $\vc f:\mathbb R^N\to\mathbb R^N$ is a known Lipschitz continuous vector field. We denote the solution of the system by $\vc u(t;\vc u_0)=\fm^t(\vc u_0)$ where $\fm^t:\mathbb R^N\to\mathbb R^N$ is the flow map associated with~\eqref{eq:ODE}. For notational simplicity, we write $\vc u(t)$ and omit its dependence on the initial condition $\vc u_0$. The corresponding observational time series are given by $\vc y(t) = S^\top \vc u(t)$. If the dynamical system is dissipative, its trajectories converge to an attractor $\mathcal A\subset \mathbb R^N$ whose dimension is smaller than $N$. If such a lower dimensional attractor does not exist, we set $\mathcal A=\mathbb R^N$.

In our setting, the initial condition $\vc u_0$ is not completely known. Instead, only its entries corresponding to the sensor locations are known which yield the observations $\vc y(0)=S^\top \vc u_0$ at the initial time. As a result, ODE~\eqref{eq:ODE} cannot be directly integrated to estimate the future states $\vc u(t)$.

Nonetheless, at any time instance $t$, S-DEIM can be used to compute the state estimation,
\begin{equation}
\tilde{\vc u}(\vc z(t)) = \Phi \vc c(\vc z(t)) = \Phi \left(S^\top \Phi\right)^+ \vc y(t) + \Phi \vc z(t),
\label{eq:S-DEIM_DA}
\end{equation}
where $\vc z(t) \in\Null[S^\top \Phi]$ for all times $t\geq 0$. The kernel vector $\vc z(t)$ is still unknown and our task is to determine its optimal value. To this end, standard data assimilation methods, such as three- and four-dimensional variational data assimilation (3DVAR and 4DVAR, respectively), can be used (cf. Ref.~\cite{Bannister2017} for a review). 

However, these data assimilation methods rely on descent iterations to minimize a suitable cost function. At each iteration, the descent direction needs to be approximated which requires repeated solves of the ODE~\eqref{eq:ODE}. In the case of 4DVAR, an additional adjoint ODE also needs to be solved at each iteration~\cite{Bannister2017}. Since the dimension $N$ of the system is often large, this leads to a high computational cost. Furthermore, the cost function is a non-convex function of the kernel vector $\vc z$; therefore, the optimization may never converge to the global minimizer. The high computational cost and non-convexity of the existing data assimilation methods hinders their applicability to high-dimensional systems. Here, we introduce an alternative approach which circumvents these two issues. 
Our method seeks $\vc z(t)$ such that the corresponding S-DEIM reconstruction $\tilde{\vc u}(\vc z(t))$ solves the ODE~\eqref{eq:ODE} as closely as possible. 

More precisely, consider a path $\vc z(t)\in\Null[S^\top\Phi]$ and its corresponding
S-DEIM reconstruction~\eqref{eq:S-DEIM_DA}. Since S-DEIM is an approximation, $\tilde{\vc u}(\vc z(t))$ does not necessarily solve the ODE~\eqref{eq:ODE}. 
Instead, we seek $\vc z(t)$ such that the error between the left- and right-hand side of the ODE is instantaneously minimized. In particular, we solve the optimization problem, 
\begin{equation}
\min_{\dot{\vc z}(t)} \left\| \frac{\id}{\id t}\tilde{\vc u}(\vc z(t)) - \vc f(\tilde{\vc u}(\vc z(t)))\right\|^2,
\label{eq:opt_inst}
\end{equation}
where
\begin{equation}
\frac{\id}{\id t}\tilde{\vc u}(\vc z(t)) = \Phi \left(S^\top \Phi\right)^+ \dot{\vc y}(t) + \Phi\dot{\vc z}(t).
\label{eq:udot_tilde}
\end{equation}

At first glance, this approach may seem problematic since the time derivative of the S-DEIM reconstruction
involves the time derivative of the observations $\dot{\vc y}$. As such, any noise in the observations can be detrimental. However, the next theorem shows that, not only the solution to the optimization problem~\eqref{eq:opt_inst} can be written explicitly, but also that the solution is independent of the derivative of the observations. 

\begin{thm}\label{thm:xiDot}
	Optimization problem~\eqref{eq:opt_inst} has a unique solution given by $\vc z(t) = Z\bxi(t)$ where $\bxi(t)\in\mathbb R^{m-n}$ solves the ODE,
	\begin{equation}
	\dot \bxi = Z^\top \Phi^\top \vc f(\tilde{\vc u}(\bxi) ),
	\label{eq:xiDot}
	\end{equation}
	where $Z$ is a kernel matrix and $\tilde{\vc u}(\bxi)$ is shorthand notation for $\tilde{\vc u}(Z\bxi)$.
\end{thm}
\begin{proof}
	Recall that the columns of $Z$ are orthonormal and $\R[Z] = \Null[S^\top\Phi]$. Consequently, for any $\vc z(t) \in\Null[S^\top\Phi]$, there exist a unique $\bxi(t)\in\mathbb R^{m-n}$ such that $\vc z(t) = Z\bxi(t)$.
	Therefore, using equation~\eqref{eq:udot_tilde}, optimization problem~\eqref{eq:opt_inst} can be written equivalently as
	\begin{equation}
	\min_{\dot{\bxi}(t)\in\mathbb R^{m-n}} \left\| \Phi \left(S^\top \Phi\right)^+ \dot{\vc y}(t) + \Phi Z\dot{\bxi}(t) - \vc f(\tilde{\vc u}(\bxi(t)))\right\|^2. 
	\end{equation}
	Since $\Phi Z$ has full column rank and $(\Phi Z)^+=Z^\top\Phi^\top$, the minimizer is unique and is given by
	\begin{equation}
	\dot \bxi = -Z^\top\Phi^\top \Phi \left(S^\top \Phi\right)^+ \dot{\vc y}(t) + Z^\top \Phi^\top \vc f(\tilde{\vc u}(\bxi) ). 
	\end{equation}
	Finally, since $\Phi^\top \Phi = \mathbb I_m$ and $Z^\top \left(S^\top \Phi\right)^+=0$ (due to Lemma~\ref{lem:orth_decomp}), we obtain the desired result.
\end{proof}

We refer to equation~\eqref{eq:xiDot} as the \emph{kernel ODE} and the corresponding approximation~\eqref{eq:S-DEIM_DA} as the DAS-DEIM reconstruction. Since the solution~\eqref{eq:xiDot} minimizes the instantaneous error between the 
left- and right-hand side of the governing ODE~\eqref{eq:ODE}, one expects that, over time, the corresponding DAS-DEIM reconstruction $\tilde{\vc u}(\bxi(t))$ evolves closer to the true state $\vc u(t)$. More precisely, consider the observations $\{\vc y(t),\allowbreak 0\leq t\leq T\}$ from the initial time to the current time $T$. Using this observations in its right-hand side, the ODE~\eqref{eq:xiDot} can be solved numerically from an initial guess $\bxi(0)=\bxi_0\in\mathbb R^{m-n}$ for the kernel vector. Even if the initial guess $\bxi_0$ is not suitable, we expect that over time it converges towards its optimal value so that the reconstruction error $\|\tilde{\vc u}(\bxi(t))-\vc u(t)\|$  decreases as time $t$ increases. This procedure is summarized in Algorithm~\ref{alg:DAS_DEIM}.
\begin{algorithm}
	\caption{Data assimilated sparse DEIM (DAS-DEIM) algorithm.}\label{alg:DAS_DEIM}
	\begin{algorithmic}
		\State{\bf Inputs:} 
		\State \textbullet~ Basis matrix $\Phi=[\bphi_1|\bphi_2|\cdots|\bphi_m]\in\mathbb R^{N\times m}$
		\State \textbullet~ Observations $\vc y(t)\in\mathbb R^n,\quad 0\leq t\leq T$
		\State \textbullet~ Vector field $\vc f:\mathbb R^N\to\mathbb R^N$
		\ \\
		\State 1. S = \texttt{pivotQR}($\Phi$) \Comment{Compute the selection matrix using~\eqref{eq:QRalg}}
		\State 2. Z = \texttt{null}($S^\top \Phi$) \Comment{Compute the kernel matrix}
		\State 3. Choose an initial guess: $\bxi_0\in\mathbb R^{m-n}$
		\State 4. Solve $\dot{\bxi} = Z^\top\Phi^\top \vc f(\tilde{\vc u}(\bxi))$ with $\bxi(0)=\bxi_0$. \Comment{Solve for $\bxi(t)$, $0\leq t\leq T$}
		\State 5. State reconstruction: $\tilde{\vc u}(t) = \Phi \left(S^\top \Phi\right)^+ \vc y(t) + \Phi Z\bxi(t),\quad 0\leq t\leq T$ 
		\ \\
		\State{\bf Outputs:} DAS-DEIM reconstruction $\tilde{\vc u}(t),\quad 0\leq t\leq T$
	\end{algorithmic}
\end{algorithm}

The dimension of the kernel ODE~\eqref{eq:xiDot} is $m-n$ which is much smaller than the dimension $N$ of the original ODE~\eqref{eq:ODE}. While standard data assimilation methods, such as 3DVAR and 4DVAR, require repeated solves of the $N$-dimensional ODE~\eqref{eq:ODE}, DAS-DEIM requires solving the smaller kernel ODE only once. As a result, DAS-DEIM is computationally more efficient than standard data assimilation techniques. Nonetheless, numerically integrating the kernel ODE~\eqref{eq:xiDot} still requires the evaluation of the $N$-dimensional vector field $\vc f$.

We show in section~\ref{sec:numerics}, with two numerical examples, that the reconstruction error of DAS-DEIM does indeed decrease over time. But before presenting our numerical results, we prove that, under certain conditions, the reconstruction error of DAS-DEIM converges exponentially fast towards zero as $t\to\infty$.

\subsection{Convergence of DAS-DEIM}\label{sec:DAS-DEIM_conv}
In order to state our main convergence result, we recall the notion of one-sided Lipschitz continuity. 
\begin{defn}\label{def:1sidedLip}
	A vector field $\vc f:\mathbb R^N\to\mathbb R^N$ is one-sided Lipschitz continuous if there exists $\rho_f\in\mathbb R$ such that
	\begin{equation}
	\langle \tilde{\vc u} -  \vc u,\vc f(\tilde{\vc u})-\vc f(\vc u)\rangle \leq \rho_f\|\tilde{\vc u}-\vc u\|^2,
	\end{equation}
	for all $\tilde{\vc u},\vc u\in\mathbb R^N$. The coefficient $\rho_f$ is referred to as the one-sided Lipschitz constant.
\end{defn}

\begin{rem}
	A few remarks about Definition~\ref{def:1sidedLip} are in order:
	\begin{enumerate}
		\item Note that $\tilde{\vc u}$ is any vector in $\mathbb R^N$ and not necessarily the S-DEIM reconstruction. Later we will apply this definition with the special case where $\tilde{\vc u}$ is the S-DEIM reconstruction of the state $\vc u$.
		\item Unlike the Lipschitz constant, the one-sided Lipschitz constant $\rho_f$ is not necessarily positive~\cite{Hu2006}. In fact, it can be zero or even a negative number. For instance, consider the one-dimensional function $f(x) = e^{-x}: (-\infty,0]\to[1,\infty)$. Since $e^{-x}\geq 1$ on its domain $(-\infty,0]$, it is easy to see that 
		$$ (x-y)(e^{-x}-e^{-y})\leq -(x-y)^2,\quad \forall x,y\in(-\infty,0]$$
		and therefore $\rho_f=-1$. If $e^{-x}$ is considered on the entire real line, then we have $\rho_f=0$.
		\item Any Lipschitz function is also one-sided Lipschitz~\cite{Hu2006}, but the converse is not necessarily the case. For instance, the exponential function $e^{-x}$ is not Lipschitz continuous, but it is one-sided Lipschitz as shown above.
	\end{enumerate}	
\end{rem}

As in Lipschitz continuous functions, the one-sided Lipschitz constant is not unique. From now on, we consider the best, i.e., the smallest possible, value of this constant. In other words, we define the one-sided Lipschitz constant as
\begin{equation}
\rho_f = \sup_{\tilde{\vc u}\neq\vc u}\frac{\langle \tilde{\vc u} -  \vc u,\vc f(\tilde{\vc u})-\vc f(\vc u)\rangle}{\|\tilde{\vc u}-\vc u\|^2}.
\end{equation}

The following result shows that, if $\vc f$ is Lipschitz continuous, both $\vc f$ and any orthogonal projection of it are one-sided Lipschitz continuous.

\begin{lem}\label{lem:1sidedLip}
	Assume $\vc f:\mathbb R^N\to\mathbb R^N$ is Lipschitz continuous with a constant $L_f>0$ and define $\vc g(\vc u)=\mathbb P\vc f(\vc u)$, for an arbitrary orthogonal projection $\mathbb P$. Then $\vc g$ is one-sided Lipschitz continuous with a constant $\rho_g\leq L_f$.
\end{lem}
\begin{proof}
	If $\mathbb P$ is the zero matrix, the result holds trivially with $\rho_g=0$. We now consider non-zero orthogonal projections.
	We have
	\begin{align}
	\rho_g = \sup_{\tilde{\vc u}\neq\vc u}\frac{\langle \tilde{\vc u} -  \vc u,\mathbb P\left[\vc f(\tilde{\vc u})-\vc f(\vc u)\right]\rangle}{\|\tilde{\vc u}-\vc u\|^2}
	&\leq \sup_{\tilde{\vc u}\neq\vc u}\frac{\|\tilde{\vc u} -  \vc u\|\cdot \|\mathbb P\|_2\cdot\|\vc f(\tilde{\vc u})-\vc f(\vc u)\|}{\|\tilde{\vc u}-\vc u\|^2}\nonumber\\
	&\leq \|\mathbb P\|_2 \sup_{\tilde{\vc u}\neq\vc u}\frac{\|\vc f(\tilde{\vc u})-\vc f(\vc u)\|}{\|\tilde{\vc u}-\vc u\|}\leq L_f,
	\end{align}
	where we used the Cauchy--Schwartz inequality and $\|\mathbb P\|_2$ denotes the spectral norm of $\mathbb P$. Since $\mathbb P$ is an orthogonal projection, we have $\|\mathbb P\|_2=1$.
\end{proof}

The next theorem is our main convergence result which uses Lemma~\ref{lem:1sidedLip} with the special choice of $\mathbb P = \Phi ZZ^\top\Phi^\top$ 
for the orthogonal projection.

\begin{thm}[Convergence of DAS-DEIM]\label{thm:DAS-DEIM_conv}
	Consider the system~\eqref{eq:ODE} and the orthogonal projection $\mathbb P=\Phi ZZ^\top\Phi^\top$. Assume that $\mathcal A\subset \R[\Phi]$ and that the one-sided Lipschitz constant $\rho_g$ of the vector field $\vc g := \mathbb P\vc f$ is negative, $\rho_g<0$.
	
	If $\vc u_0\in\mathcal A$, then 
	\begin{equation}
	\lim_{t\to\infty}\|\tilde{\vc u}(t)- \vc u(t)\|=0,
	\end{equation}
	where $\tilde{\vc u}(t)=\Phi (S^\top\Phi)^+\vc y(t) + \Phi Z\bxi(t)$ is the S-DEIM reconstruction with $\bxi(t)$ satisfying equation~\eqref{eq:xiDot}. Furthermore, the convergence is exponentially fast.
\end{thm}
\begin{proof}
	Since $\vc u_0\in\mathcal A$ and the attractor is invariant, we have $\vc u(t)\in\mathcal A$ for all $t\geq 0$. Since $\mathcal A\subset \R[\Phi]$, using Corollary~\ref{cor:exactReconst}, we have $\vc u(t)\in\R[\Phi]$ and therefore, 
	\begin{equation}
	\vc u(t) = \Phi(S^\top\Phi)^+\vc y(t) + \Phi ZZ^\top\Phi^\top \vc u(t),
	\end{equation}
	for all $t\geq 0$. Taking the time derivative of this equation, we obtain
	\begin{equation}\label{eq:udot}
	\dot{\vc u} = \Phi(S^\top\Phi)^+\dot{\vc y} + \Phi ZZ^\top\Phi^\top \vc f(\vc u).
	\end{equation}
	
	Now let us consider the time derivative of the S-DEIM reconstruction $\tilde{\vc u}$,
	\begin{equation}\label{eq:udot_tilde2}
	\dot{\tilde{\vc u}} = \Phi(S^\top\Phi)^+\dot{\vc y} + \Phi ZZ^\top\Phi^\top \vc f(\tilde{\vc u}),
	\end{equation}
	where we used the fact that $\bxi(t)$ satisfies equation~\eqref{eq:xiDot}.
	Subtracting equations~\eqref{eq:udot} and~\eqref{eq:udot_tilde2}, we obtain
	\begin{equation}
	\frac{\id}{\id t}\left(\tilde{\vc u} - \vc u\right) = \mathbb P \left[ \vc f(\tilde{\vc u})-\vc f(\vc u)\right],
	\end{equation}
	where $\mathbb P =\Phi ZZ^\top\Phi^\top$ is an orthogonal projection. Taking the inner product with $\tilde{\vc u} - \vc u$, we obtain
	\begin{equation}
	\frac12 \frac{\id}{\id t}\|\tilde{\vc u} - \vc u\|^2 = \langle \tilde{\vc u} - \vc u,\mathbb P\left[ \vc f(\tilde{\vc u})-\vc f(\vc u)\right]\rangle
	\leq \rho_g \|\tilde{\vc u} - \vc u\|^2,
	\end{equation}
	where we used the one-sided Lipschitz continuity of $\vc g = \mathbb P\vc f$. 
	
	Finally, using the Gr\"onwall's inequality, we have
	\begin{equation}
	\|\tilde{\vc u}(t) - \vc u(t)\|\leq \|\tilde{\vc u}(0) - \vc u(0)\|e^{\rho_g t}.
	\end{equation}
	Since the one-sided Lipschitz constant $\rho_g$ is negative, we obtain the desired result that the error tends to zero as $t\to\infty$. Clearly, this convergence is exponential in time.
\end{proof}

We conclude this section with a few remarks regarding Theorem~\ref{thm:DAS-DEIM_conv}. As in case of Lipschitz constant, estimating the smallest one-sided Lipschitz constant is difficult in practice. We also note that the conditions of this theorem are sufficient but not necessary. As our numerical results in section~\ref{sec:numerics} indicate, DAS-DEIM converges for complex dynamical systems, where these conditions are not necessarily met. 

Finally, we point out that Theorem~\ref{thm:DAS-DEIM_conv} can be used to inform the choice of the basis matrix $\Phi$ and the selection matrix $S$. Note that the orthogonal projection  $\mathbb P=\Phi ZZ^\top\Phi^\top$ depends on these matrices. Therefore, in principle, it might be possible to choose these matrices in such a way that the vector field $\vc g = \mathbb P\vc f$ has a negative one-sided Lipschitz constant.

\section{Numerical results}\label{sec:numerics}
We present two numerical examples using Lorenz63 and Lorenz96 systems. For each system, two trajectories are numerically generated from different initial conditions. The first trajectory is used as training data to extract POD modes $\{\bphi_i\}_{i=1}^m$ and to compute the selection matrix $S$ using the QR algorithm~\eqref{eq:QRalg}. The second trajectory is used for testing purposes, i.e., gather observational data, solve the kernel ODE~\eqref{eq:xiDot}, and form the corresponding DAS-DEIM estimation~\eqref{eq:S-DEIM_DA}. All errors reported in this section correspond to the testing data set. The observations are recorded at time intervals of $\Delta t=0.2$. To solve the kernel ODE~\eqref{eq:xiDot}, piecewise linear interpolation is used to approximate the observations in intermediate time instances. In the DAS-DEIM Algorithm~\ref{alg:DAS_DEIM}, we use the initial guess $\bxi_0=\vc 0$.

We reiterate that (Q)-DEIM was originally developed for reduced-order modeling~\cite{Sorensen2010,Drmac2016}. Here, however, our goal is estimating the full state of the system from its partial sparse observations without any model order reduction.

\subsection{Lorenz63 system}
\begin{figure}
	\centering
	\includegraphics[width=\textwidth]{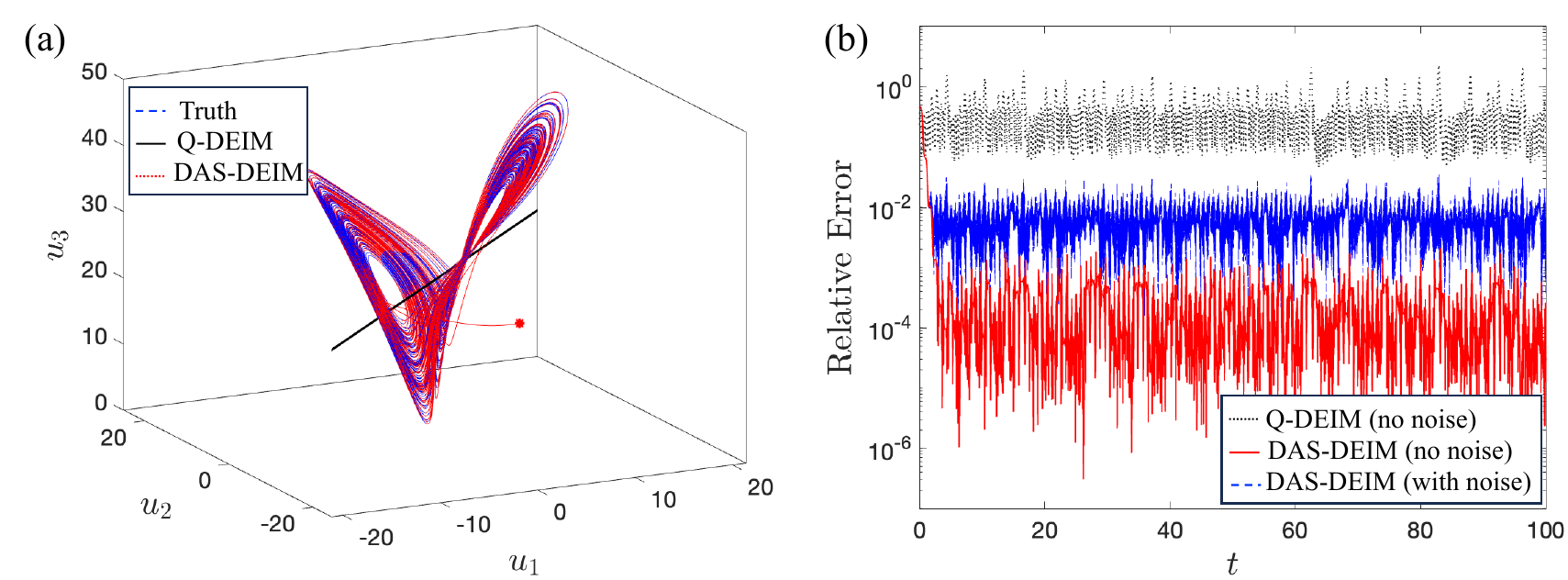}
	\caption{Vanilla DEIM and S-DEIM state estimation for Lorenz63. (a) State space showing the true trajectory (blue), vanilla Q-DEIM (straight black line), and DAS-DEIM (red). The red circle marks the initial  DAS-DEIM estimation. (b) Relative errors as a function of time.}
	\label{fig:Lorenz63}
\end{figure}

As a benchmark example, we first consider the Lorenz attractor (or Lorenz63 system), 
\begin{equation}
\dot u_1 = \sigma(u_2-u_1),\quad \dot u_2 = u_1(\rho-u_3)-u_2,\quad \dot u_3=u_1u_2-\beta u_3,
\end{equation}
with the parameters $\sigma = 10$, $\rho=28$, and $\beta=8/3$ which is known to exhibit chaotic dynamics. We assume that only $n=1$ observation is available. The QR algorithm~\eqref{eq:QRalg} determines the second component of the system as most informative; therefore, we set $y(t) = u_2(t)$. We first apply vanilla DEIM with $m=n=1$ POD mode. As shown in Figure~\ref{fig:Lorenz63}(a), since the range of $\Phi$ is one-dimensional, the estimated state evolves on a line segment. As a result, the relative error is quite large with a mean of approximately 36\%.

For S-DEIM, we can increase the number of modes arbitrarily. Therefore, we use $m=3$ modes with the same observational time series as before ($y(t) = u_2(t)$). As shown in figure~\ref{fig:Lorenz63}(b), the initial relative error for DAS-DEIM is similar to vanilla DEIM. However, the DAS-DEIM error decreases rapidly as the kernel vector $\vc z(t)$ converges towards its optimum. Eventually, this error oscillates around $10^{-4}$ (0.01\%), indicating excellent agreement with the true solution of the Lorenz63 system.

We also experiment with noisy data. In particular, we repeat the DAS-DEIM estimation by adding an error $\varepsilon(t)\sim N(0,0.1)$ to the observations $y(t)=u_2(t)$ at each time instance. Here, $N(0,0.1)$ denotes the normal distribution with mean zero and standard deviation $0.1$. As shown in figure~\ref{fig:Lorenz63}(b), the noise negatively affects the estimation so that the relative error increases, with a mean of 0.7\%. Nonetheless, this error is still quite small and significantly smaller than vanilla DEIM error of 36\% with clean data.

One may wonder what happens if we use vanilla DEIM with $n=1$ observations (as before) and $m=3$ modes. In this case, the vanilla DEIM estimation error increases to 43\% (not shown here). This is in line with the numerical results of~\cite{Brunton2021} who observed that the minimum vanilla DEIM error occurs when the number of modes is smaller than or equal to the number of sensors ($m\leq n$).

\subsection{Lorenz96 system}
\begin{figure}
	\centering
	\includegraphics[width=\textwidth]{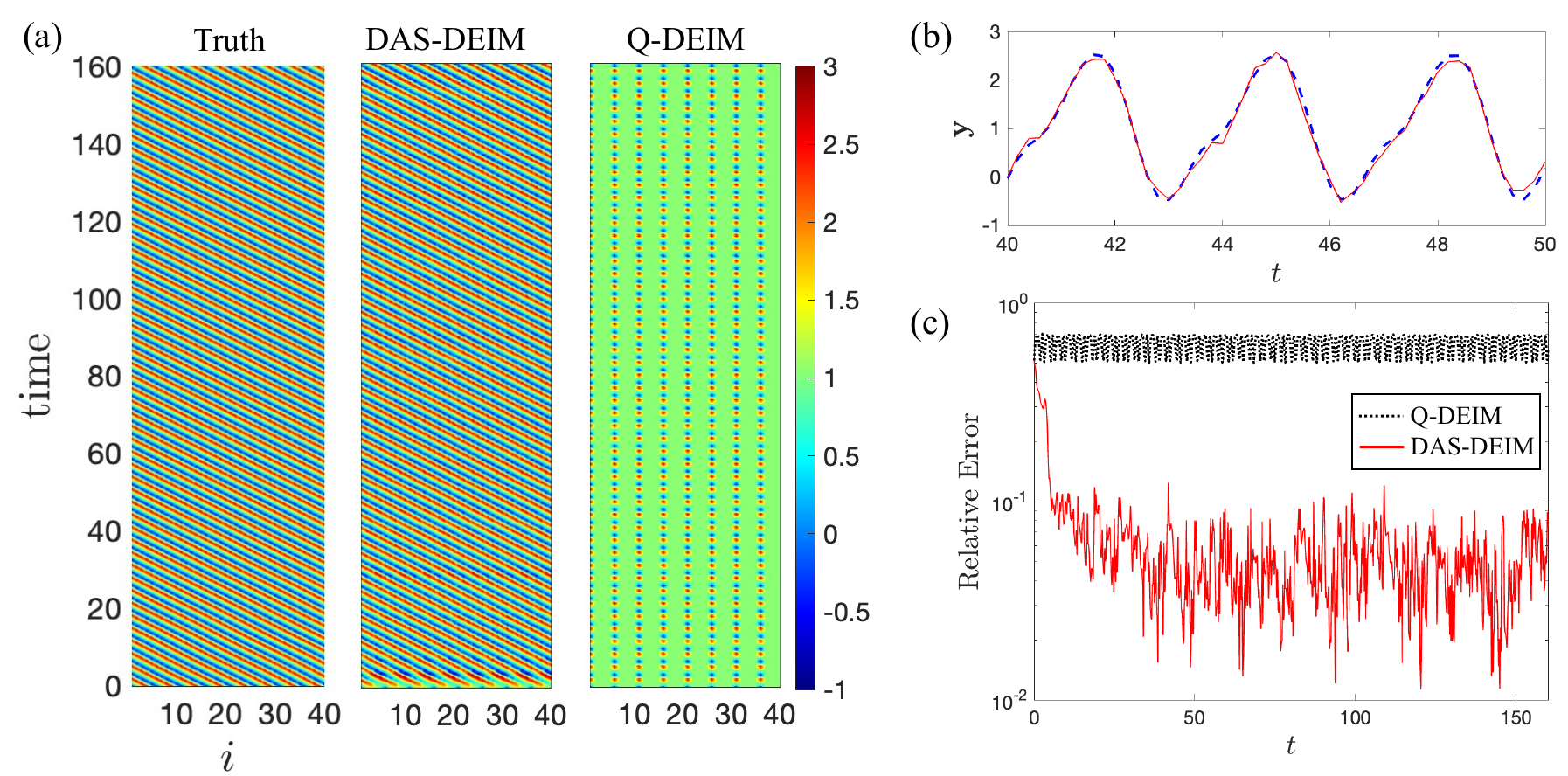}
	\caption{State estimation for Lorenz96 system. (a) True solution, DAS-DEIM, and vanilla Q-DEIM reconstructions. (b) Close-up view of the observational data from $n=1$ sensor. Dashed blue shows the true observation data and solid red line marks the corresponding noisy observations. (c) Relative error of vanilla Q-DEIM compared to DAS-DEIM.}
	\label{fig:Lorenz96}
\end{figure}

In this section, we consider the Lorenz96 system,
\begin{equation}
\dot u_i = (u_{i+1}-u_{i-2})u_{i-1}-u_i+F,\quad 1\leq i\leq N,
\end{equation}
with $N=40$ and periodic boundary conditions: $u_{-1} = u_{N-1}$, $u_0=u_N$, $u_{N+1}=u_1$. This system is a conceptual model for atmospheric flow and is routinely used for assessing data assimilation methods~\cite{Chattopadhyay2020,Lguensat2017}. Here we use the forcing amplitude of $F=2$ which leads to a modulating traveling wave solution as shown in figure~\ref{fig:Lorenz96}(a).

The singular values corresponding to the POD modes exhibit a sharp decay after the first five modes: the fifth singular value is around 0.9 whereas the sixth singular value is approximately $4\times 10^{-8}$. This indicates that the attractor $\mathcal A$ of the system resides approximately in the linear subspace spanned by the first five POD modes.
Therefore, we use $m=5$ POD modes for both vanilla DEIM and DAS-DEIM estimations. In order to demonstrate that DAS-DEIM performs well even if very few sensors are available, we use $n=1$ sensor measurements. Column pivoted QR factorization algorithm~\eqref{eq:QRalg} determines $u_1$ as the best place for sensor placement. 
Figure~\ref{fig:Lorenz96}(b) shows a close-up view of the observational data $y(t)=u_1(t)+\varepsilon(t)$ where random noise $\varepsilon\sim N(0,0.1)$ is added to the true data at each time instance. Both vanilla DEIM and DAS-DEIM use this noisy data as input.

As shown in figure~\ref{fig:Lorenz96}(a), vanilla DEIM reproduces incorrect dynamics. In contrast, after initial transients, DAS-DEIM reconstruction converges towards the true solution of the system. Figure~\ref{fig:Lorenz96}(c) shows that the vanilla DEIM error oscillates around 62\%, whereas DAS-DEIM error decays rapidly and eventually oscillates around 5\%. Recall that DAS-DEIM requires the solution of the kernel ODE~\eqref{eq:xiDot} whose dimension $m-n=4$ is much smaller than the dimension $N=40$ of the Lorenz96 system. Solving this kernel ODE using Matlab's ODE45 took only 0.4 seconds.
\begin{figure}
	\centering
	\includegraphics[width=.5\textwidth]{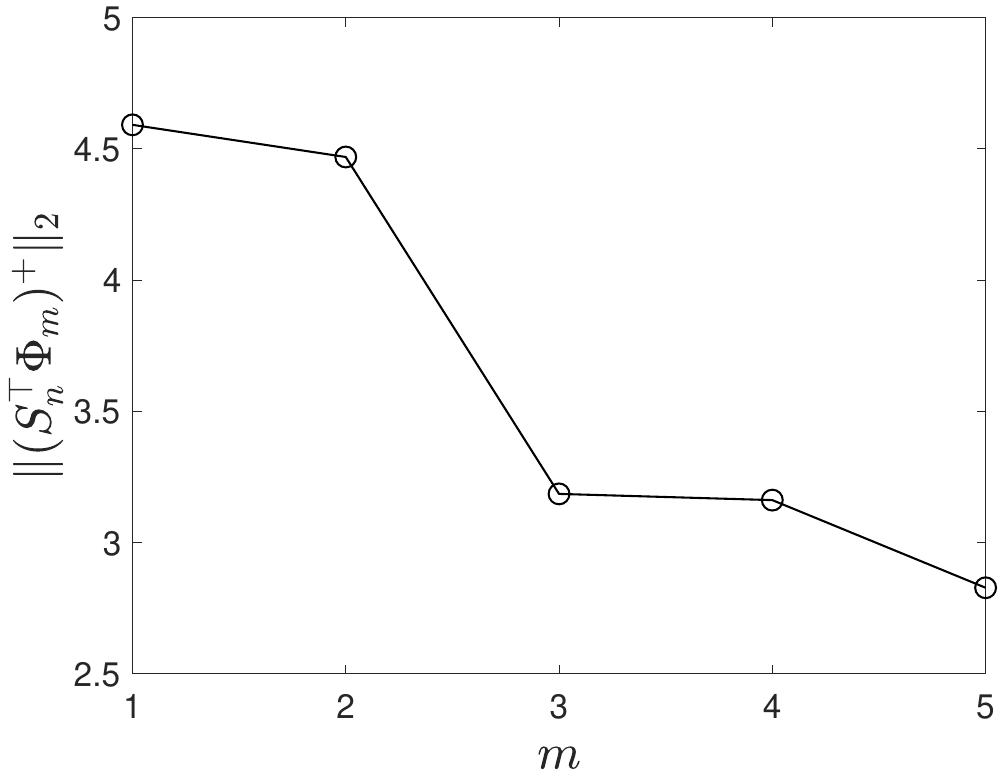}
	\caption{The matrix norm $\|(S_n^\top\Phi_m)^+\|_2$ for the Lorenz96 system with $n=1$ sensors and increasing number of modes $m$.}
	\label{fig:Lorenz96_conjecture}
\end{figure}

Finally, we recall Remark~\ref{rem:prefactor} which states that the prefactor $\|(S_n^\top\Phi_m)^+\|_2$ in the error upper bound~\eqref{eq:SDEIM_ub} is a decreasing function of $m$ (the number of modes). Figure~\ref{fig:Lorenz96_conjecture} demonstrates this numerically with $n=1$ sensor. As the number of modes $m$ increases, the prefactor $\|(S_n^\top\Phi_m)^+\|_2$ decreases indicating that the S-DEIM reconstruction should improve as more modes are used.

\section{Conclusions}\label{sec:conc}
We introduced Sparse DEIM (S-DEIM) for interpolating unknown functions from their sparse observations when the number of observations is limited. In the special case of continuous-time dynamical systems, we introduced the data assimilated S-DEIM (DAS-DEIM) algorithm to efficiently estimate the unknown kernel vector in S-DEIM. Our numerical and theoretical results show great promise for fast and accurate estimation of the state of the system from very few observational data.

A number of open problems remain to be addressed. As we mentioned at the end of section~\ref{sec:DAS-DEIM_conv}, Theorem~\ref{thm:DAS-DEIM_conv} can be used to inform sensor placement. Although here we used the standard column pivoted QR algorithm, a tailor-made sensor placement algorithm for DAS-DEIM should be investigated.
Furthermore, Theorem~\ref{thm:DAS-DEIM_conv} provides sufficient conditions for the convergence of the DAS-DEIM algorithm. These conditions are most likely too restrictive and therefore a convergence theorem with necessary and sufficient conditions is desirable.

Finally, we point out that certain data sets, such as those encountered in image processing, are not necessarily generated by a dynamical system. In such cases, alternative methods for approximating the optimal kernel vector of S-DEIM are needed. For instance, deep learning methods can potentially be used. However, our preliminary results (not shown here) indicate that simple neural network architectures (e.g., feed forward) fail to yield an acceptable approximation. In the future, the use of more advanced architectures, such as long short-term memory (LSTM) networks and autoencoders, will be explored.

\section*{Acknowledgments}
I am grateful to Profs. Ilse Ipsen and Arvind Saibaba for fruitful conversations.


\end{document}